\documentclass[journal]{IEEEtran}
\ifCLASSINFOpdf
\else
\fi

\usepackage{float}
\usepackage{caption}
\usepackage{subcaption}
\usepackage{tikz,pgfplots,pgfplotstable,booktabs}
\usepackage[official]{eurosym}
\usepackage{amsmath}
\usepackage{amssymb}
\allowdisplaybreaks[1]

\pgfplotsset{
    discard if not/.style 2 args={
        x filter/.code={
            \edef\tempa{\thisrow{#1}}
            \edef\tempb{#2}
            \ifx\tempa\tempb
            \else
                
            \fi
        }}}

\begin{document}
\title{Solving Linear Bilevel Problems Using Big-Ms: Not All That Glitters Is Gold}

\author{Salvador Pineda and
        Juan Miguel Morales
\thanks{S. Pineda is with the Department
of Electrical Engineering, University of Malaga, Spain. E-mail: spinedamorente@gmail.com.}% <-this % stops a space
\thanks{J. M. Morales is with the Department of Applied Mathematics, University of Malaga, Spain. E-mails: juan.morales@uma.es; juanmi82mg@gmail.com.}% <-this % stops a space
\thanks{This work was supported in part by the Spanish Ministry of Economy, Industry and Competitiveness through projects ENE2016-80638-R and ENE2017-83775-P. This project has received funding from the European Research Council (ERC) under the European Union’s Horizon 2020 research and innovation programme (grant agreement No 755705).}}

\maketitle

\begin{abstract}
The most common procedure to solve a linear bilevel problem in the PES community is, by far, to transform it into an equivalent single-level problem by replacing the lower level with its KKT optimality conditions. Then, the complementarity conditions are reformulated using additional binary variables and large enough constants (big-Ms) to cast the single-level problem as a mixed-integer linear program that can be solved using optimization software. In most cases, such large constants are tuned by trial and error. We show, through a counterexample, that this widely used trial-and-error approach may lead to highly suboptimal solutions. Then, further research is required to properly select big-M values to solve linear bilevel problems.
\end{abstract}


% \fi
\IEEEpeerreviewmaketitle

\section{Introduction}\label{sec:intro}

\IEEEPARstart{D}{ecentralized} environments are characterized by multiple decisions makers with divergent objectives that interact with each other. One of the simplest instances only considers two decision makers that make their decisions in a sequential manner. The player deciding first is called the \textit{leader}, while the one deciding afterwards is called the \textit{follower}. This non-cooperative sequential game is known as a \textit{Stackelberg game} and was first investigated in \cite{VonStackelberg1952}. A Stackelberg game can be mathematically formulated as a bilevel problem (BP) \cite{Bard1998,Dempe2002}. If the objective functions of both players and all constraints are linear, the resulting linear bilevel problem (LBP) can be generally formulated as follows:
\begin{subequations}
\begin{align}
 \min_{x\in \mathbb{R}^n}  \quad & a^T x+ b^T y \label{bp1_1} \\
 \text{s.t.} \quad &  c^T_i x+ d^T_i y\leq e_i \quad \forall i \label{bp1_2}  \\
& \min_{y\in \mathbb{R}^m} \quad  p^T x+ q^T y  \label{bp1_3} \\
& \,\, \text{s.t.} \quad \, r^T_j x+ s^T_j y\leq t_j \; (\lambda_j) \quad \forall j \label{bp1_4}
\end{align} \label{bp1}
\end{subequations} 
\noindent where $a,b,c_i,d_i,p,q,r_j,s_j$ and $e_i,t_j$ are vectors of appropriate dimension and scalars, respectively. The dual variable of the lower-level constraint \label{bp1_4} is denoted by $\lambda_j$ in brackets. Mathematically, upper-level constraints that include upper- and lower-level variables can lead to disconnected feasible regions \cite{Colson2005a}, which  complicates the solution of the LBP as illustrated in \cite{Shi2005a}. Dealing with the solution to this variant goes beyond the purposes of this letter and thus, we assume $d_i=0$. \textcolor{black}{This assumption is common in several applications of linear bilevel programming in the PES technical literature. For example, in long-term planning models formulated as bilevel problems \cite{Wogrin2011a,Kazempour2011c,Baringo2014a,Maurovich-Horvat2014}, the upper-level problem determines investment decisions to maximize investor's profit, while the lower-level problem yields the dispatch quantities to minimize operating cost. In most cases, upper-level constraints model maximum available capacities to be installed and/or budget limitations, but do not include lower-level dispatch variables.}

Since the lower-level optimization problem is linear, it can be replaced with its KKT optimality conditions as follows:
\begin{subequations}
\begin{align}
 \min_{x\in \mathbb{R}^n,y\in\mathbb{R}^m}  \quad & a^T x+ b^T y \label{bp2_1} \\
 \text{s.t.} \quad &  c^T_i x+ d^T_i y\leq e_i \quad \forall i  \label{bp2_2}  \\
& r^T_j x+ s^T_j y\leq t_j \quad \forall j \label{bp2_3} \\
& q + \sum_j \lambda_j s_j = 0 \label{bp2_4} \\
& \lambda_j \geq 0, \quad \forall j \label{bp2_5} \\
& \lambda_j \left( r^T_j x+ s^T_j y - t_j \right) = 0, \quad \forall j \label{bp2_6} 
\end{align} \label{bp2}
\end{subequations} 
Non-linear complementarity constraints \eqref{bp2_6} are further handled using the Fortuny-Amat mixed-integer reformulation \cite{Fortuny-Amat1981} as presented below:
\begin{subequations}
\begin{align}
 \min_{x\in \mathbb{R}^n,y\in\mathbb{R}^m}  \quad & a^T x+ b^T y \label{bp3_1} \\
 \text{s.t.} \quad & \eqref{bp2_2}-\eqref{bp2_5} \\
& \lambda_j \leq u_j M^D_j, \quad \forall j \label{bp3_5} \\ 
& -r^T_j x- s^T_j y + t_j \leq (1-u_j) M^P_j, \forall j \label{bp3_6} \\
& u_j \in \{0,1\}, \quad \forall j \label{bp3_7} 
\end{align} \label{bp3}
\end{subequations} 
\noindent where $u_j$ are additional binary variables and $M^P_j,M^D_j$ are large enough constants. Model \eqref{bp3} is a mixed-integer linear optimization problem that can be solved using commercial software. \textcolor{black}{However, formulation \eqref{bp3} is equivalent to formulation \eqref{bp2} provided that the large enough constants $M^P_j,M^D_j$ are valid upper bounds for the primal and dual variables of the lower-level problem, respectively.} Notice that appropriate values for $M^P_j$ are often available, because they relate to primal variables, which are typically bounded by nature. However, $M^D_j$ are upper bounds on dual variables and therefore, tuning these large enough constants is a more challenging task. The most commonly used trial-and-error tuning procedure reported in the technical literature runs as follows:
\begin{enumerate}
    \item Select initial values for $M^P_j$ and $M^D_j$.
    \item Solve model \eqref{bp3}.
    \item Find a $j'$ such that $u_{j'}=0$ and $-r^T_{j'} x- s^T_{j'} y + t_{j'} = M^P_{j'}$. If such a $j'$ exists, increase the value of $M^P_{j'}$ and go to step 2). Otherwise, go to step 4).
    \item Find a $j'$ such that $u_{j'}=1$ and $\lambda_{j'}=M^D_{j'}$. If such a $j'$ exists, increase the value of $M^D_{j'}$ and go to step 2). Else, the solution to \eqref{bp3} \emph{is assumed} to correspond to the optimal solution of the original bilevel problem \eqref{bp1}.
\end{enumerate}

The trial-and-error procedure described above has been used in a great number of research works in the PES technical literature related to electricity grid security analysis \cite{Motto2005}, transmission expansion planning \cite{Garces2009,Jenabi2013b}, strategic bidding of power producers \cite{Ruiz2009,Zugno2013}, generation capacity expansion \cite{Wogrin2011a,Kazempour2011c}, investment in wind power generation \cite{Baringo2014a,Maurovich-Horvat2014} and market equilibria models \cite{Pozo2011,Ruiz2012}, among many others. Furthermore, its use is likely to continue in the future. However, as shown in the next section with a simple counterexample, this trial-and-error procedure does not guarantee global optimality of the original bilevel problem.

\section{Counterexample}

Let us consider the following linear bilevel problem:
\begin{subequations}
\begin{align}
 \max_{x\in \mathbb{R}}  \quad & z = x+y \label{bp4_1} \\
 \text{s.t.} \quad &  0 \leq x \leq 2 \label{bp4_2}  \\
& \min_{y\in \mathbb{R}} \quad  y  \label{bp4_3} \\
& \,\, \text{s.t.} \quad \, y \geq 0 \quad (\lambda_1) \label{bp4_4} \\
& \quad \quad \;\; x - 0.01y \leq 1 \quad  (\lambda_2) \label{bp4_5}
\end{align} \label{bp4}
\end{subequations} 
It is easy to verify that the optimal solution to this problem is $z^*=102, x^*=2, y^*=100, \lambda_1^*=0, \lambda_2^*=100$. Following the procedure described in Section \ref{sec:intro}, we can reformulate \eqref{bp4} as the following mixed-integer linear programming problem:
\begin{subequations}
\begin{align}
 \max_{x\in \mathbb{R},y\in\mathbb{R}}  \quad & z = x+y \label{bp5_1} \\
 \text{s.t.} \quad &  0 \leq x \leq 2 \label{bp5_2}  \\
& y \geq 0 \label{bp5_4} \\
& x - 0.01y \leq 1 \label{bp5_5} \\
& 1 - \lambda_1 - 0.01\lambda_2 = 0 \label{bp5_6} \\
& \lambda_1, \lambda_2 \geq 0 \\
& \lambda_1 \leq u_1 M^D_1 \\
& y \leq (1-u_1) M^P_1 \\
& \lambda_2 \leq u_2 M^D_2 \\
& -x + 0.01y +1 \leq (1-u_2) M^P_2 \\
& u_1, u_2 \in \{0,1\}
\end{align} \label{bp5}
\end{subequations} 
\textcolor{black}{To solve problem \eqref{bp5} we follow the steps of the trial-and-error procedure described in Section \ref{sec:intro}:
\begin{enumerate}
    \item Select initial values for the large enough constants, for example, $M^P_1 = M^P_2 = 200$ and $M^D_1 = M^D_2 = 50$.
    \item Solve problem \eqref{bp5} using a mixed-integer optimization software. Alternatively, and due to the small size of our example, we can also solve problem \eqref{bp5} by simple enumeration, that is, by exploring the solutions to the four linear optimization problems that arise from all the possible 0-1 combinations of the binary variables $u_1$ and $u_2$. The results of solving these four linear programming problems are collated in Table \ref{tab:sol}, which includes the values of $x,y,\lambda_1,\lambda_2, z$.
    \begin{table}[H]
    \centering
    \begin{tabular}{|c|c|c|c|c|c|c|c|}
    \hline
        Case & $u_1$ & $u_2$ & $x$ & $y$ & $\lambda_1$ & $\lambda_2$ & $z$ 
        \\
        \hline
        1 & 0 & 1 & \multicolumn{5}{c|}{Infeasible} \\
        \hline
        2 & 1 & 1 & 1 & 0 & \multicolumn{2}{c|}{Multiple} & 1 \\
        \hline
        3 & 1 & 0 & 1 & 0 & 1 & 0 & 1 \\
        \hline
        4 & 0 & 0 & \multicolumn{5}{c|}{Infeasible} \\
    \hline
    \end{tabular}
    \caption{Results for $M^P_1 = M^P_2 = 200$, $M^D_1 = M^D_2 = 50$}
    \label{tab:sol}
\end{table}
    If $u_1$ and $u_2$ are respectively fixed to 0 and 1 (case 1), the resulting linear problem is infeasible. The same happens in case 4, where both binary variables are fixed to 0. Case 2 has a unique solution for the primal variables but multiple dual solutions that satisfy $0\leq\lambda_1\leq 50$, $0\leq\lambda_2\leq 50$ and $1-\lambda_1-0.01\lambda_2 = 0$. Finally, the solution to case 3 is a singleton, with the same objective function as that of case 2. Therefore, the solution to \eqref{bp5} corresponds to that of case 3 or of case 2, i.e., $x=1,y=0$, and $z=1$. \\
    We also solved problem \eqref{bp5} using \textcolor{black}{CPLEX 12.7.0.0} under GAMS. We tested above 20 combinations of CPLEX options regarding the method used to solve such a mixed-integer linear program and in all cases the obtained solution was that of case 3, namely, $x=1, y=0, \lambda_1=1, \lambda_2=0$, and $z=1$.
    \item Since $u_2=0$ but $-1 + 0.01\cdot 0 +1 < 200$, then $M^P_2$ does not need to be increased and we go to step 4).
    \item Since $u_1=1$ but $1<50$, then $M^D_1$ does not need to be increased and the solution found in step 2) is assumed to be the optimal solution to the original linear bilevel problem \eqref{bp4}.
\end{enumerate}}

\textcolor{black}{
Notice that the solution provided \emph{as optimal} by the trial-and-error procedure ($z=1$) is far lower than the actual optimal one ($z^*=102$). This happens because the  large constant $M^D_2$ is arbitrarily set to 50, which is lower than the actual optimal value of the corresponding dual variable $\lambda_2$, which is 100. And what is more relevant, the trial-and-error procedure is unable to detect such a poor tuning of the Big-Ms.
}

\textcolor{black}{
Next, let us analyze the results of the trial-and-error procedure if the large enough constants are somehow set to appropriate values:
\begin{enumerate}
 \item Select initial values for the large enough constants, for instance, $M^P_1 = M^P_2 = 200$ and $M^D_1 = M^D_2 = 200$.
 \item Solve problem \eqref{bp5}. Similarly to Table \ref{tab:sol}, Table \ref{tab:sol2} provides the values taken on by the objective function $z$ and the primal and dual variables of problem \eqref{bp5} for every possible combination of the binary variables $u_1$ and $u_2$. 
    \begin{table}[H]
    \centering
    \begin{tabular}{|c|c|c|c|c|c|c|c|}
    \hline
        Case & $u_1$ & $u_2$ & $x$ & $y$ & $\lambda_1$ & $\lambda_2$ & $z$ 
        \\
        \hline
        1 & 0 & 1 & 2 & 100 & 0 & 100 & 102 \\
        \hline
        2 & 1 & 1 & 1 & 0 & \multicolumn{2}{c|}{Multiple} & 1 \\
        \hline
        3 & 1 & 0 & 1 & 0 & 1 & 0 & 1 \\
        \hline
        4 & 0 & 0 & \multicolumn{5}{c|}{Infeasible} \\
    \hline
    \end{tabular}
\caption{Results for $M^P_1 = M^P_2 = M^D_1 = M^D_2 = 200$}
    \label{tab:sol2}
\end{table}
    The optimal solution is then $x=2, y=100, \lambda_1=0, \lambda_2=100, z=102$.
    \item Since $u_1=0$ but $100 < 200$, then $M^P_1$ does not need to be increased and we go to step 4).
    \item Since $u_2=1$ but $100<200$, then $M^D_2$ does not need to be increased and the solution found in step 2) is the optimal solution to problem \eqref{bp4}.
\end{enumerate}
}

\textcolor{black}{
This counterexample confirms that the trial-and-error procedure used in the technical literature and, in particular, within the PES community, only works if the big-Ms are, in some way, set to appropriate values. Otherwise, the trial-and-error procedure may fail and provide highly suboptimal solutions. Notice, however, that it is particularly challenging to guess proper bounds for the dual variables of the lower-level problem. Therefore, caution must be exercised when using such a trial-and-error procedure. Most importantly, one must bear in mind that the solution this procedure provides cannot be guaranteed to be the global optimal solution to the original linear bilevel problem.
}

\textcolor{black}{The ultimate aim of this counterexample is simply to illustrate the drawbacks of this commonly used trial-and-error procedure. Notwithstanding this, similar suboptimal results can also be obtained in more realistic applications of linear bilevel problems that have been considered in the technical literature. In fact, we have experienced this issue in generic linear bilevel problems of higher dimension for which the lower-level dual variables take on values of a very different order of magnitude since, in such cases, properly tuning the large constants becomes particularly difficult.}

\textcolor{black}{One possibility to increase the chances of success of the described trial-and-error method is to better estimate valid bounds for the lower-level dual variables. In this line, authors of \cite{Pineda2017} propose a method that first solves problem \eqref{bp2} as a non-linear optimization problem that includes neither binary variables nor large constants. Although the obtained solution is only locally optimal, the order of magnitude of primal and dual variables for such a local optimal solution can be used to subsequently adjust the values of the large constants in problem \eqref{bp3}. Results presented in \cite{Pineda2017} for a large set of randomly generated examples show the good performance of such a method to select big-M values.}

\section{Conclusion}

This letter aims to raise concern about the widespread and continued use of the big-M approach to solve LBP within the PES community. We show, using a counterexample, that the trial-and-error procedure that is presently employed to tune the big Ms in many works published in PES journals may actually fail and provide highly suboptimal solutions. We advocate, instead, for the use of more sophisticated methods like the one proposed in \cite{Pineda2017} to properly tune the values of the big-Ms when solving LBP.

% or


%\appendices

%\section{}

%\section*{Acknowledgment}

%This work was supported in part by the Spanish Ministry of Economy, Industry and Competitiveness through project ENE2016-80638-R and in part by the Research Funding Program for Young Talented Researchers of the University of M\'{a}laga through project PPIT-UMA-B1-2017/18.

\ifCLASSOPTIONcaptionsoff
  \newpage
\fi



\bibliographystyle{IEEEtran}
\bibliography{mendeley}
\end{document}